\documentclass{amsart}
\usepackage[active]{srcltx}
\usepackage{amsthm,mathrsfs}
\usepackage{a4wide}
\usepackage[english]{babel}
\usepackage{amsfonts}
\usepackage{amsmath}
\usepackage{amssymb}
\usepackage{dsfont}

\usepackage{float}
\usepackage{url}

\makeatletter\def\blfootnote{\xdef\@thefnmark{}\@footnotetext}\makeatother
\allowdisplaybreaks
\title{\bf Metric number theory, lacunary series and systems
of dilated functions}

\author{Christoph Aistleitner} 
\address{Department of Applied Mathematics, School of Mathematics and Statistics, University of New South Wales, Sydney NSW 2052, Australia}
\email{aistleitner@math.tugraz.at}

\thanks{The author is supported by a Schr\"odinger scholarship of the Austrian Research
Foundation (FWF)}


\begin{document}

\begin{abstract}
By a classical result of Weyl, for any increasing sequence $(n_k)_{k \geq 1}$ of integers the sequence of fractional parts $(\{n_k x\})_{k \geq 1}$ is uniformly distributed modulo 1 for almost all $x \in [0,1]$. Except for a few special cases, e.g. when $n_k=k,~k \geq 1$, the exceptional set cannot be described explicitly. The exact asymptotic order of the discrepancy of $(\{n_k x\})_{k \geq 1}$ for almost all $x$ is only known in a few special cases, for example when $(n_k)_{k \geq 1}$ is a (Hadamard) lacunary sequence, that is when $n_{k+1}/n_k \geq q > 1,~k \geq 1$. In this case of quickly increasing $(n_k)_{k \geq 1}$ the system $(\{n_k x\})_{k \geq 1}$ (or, more generally, $(f(n_k x))_{k \geq 1}$ for a 1-periodic function $f$) shows many asymptotic properties which are typical for the behavior of systems of \emph{independent} random variables. Precise results depend on a fascinating interplay between analytic, probabilistic and number-theoretic phenomena.\\

Without any growth conditions on $(n_k)_{k \geq 1}$ the situation becomes much more complicated, and the system $(f(n_k x))_{k \geq 1}$ will typically fail to satisfy probabilistic limit theorems. An important problem which remains is to study the almost everywhere convergence of series $\sum_{k=1}^\infty c_k f(k x)$, which is closely related to finding upper bounds for maximal $L^2$-norms of the form
$$
\int_0^1 \left( \max_{1 \leq M \leq N} \left| \sum_{k=1}^M c_k f(kx) \right| \right)^2 dx.
$$ 
The most striking example of this connection is the equivalence of the Carleson convergence theorem and the Carleson--Hunt inequality for maximal partial sums of Fourier series. For general functions $f$ this is a very difficult problem, which is related to finding upper bounds for certain sums involving greatest common divisors.
\end{abstract}

\date{}
\maketitle

\section{Uniform distribution modulo 1} \label{sec1}

A sequence of real numbers from the unit interval $(x_k)_{k \geq 1}$ is called \emph{uniformly distributed modulo 1} (u.d. mod 1) if for all $0 \leq a < b \leq 1$ the asymptotic relation 
\begin{equation} \label{ud}
\lim_{N \to \infty} \frac{1}{N} \sum_{k=1}^N \mathbf{1}_{[a,b)} (x_k) = b-a
\end{equation}
holds. Roughly speaking, a sequence is u.d. mod 1 if asymptotically every interval $[a,b) \subset [0,1]$ receives its fair share of points, which is proportional to its length. In an informal way, uniformly distributed sequences are often considered as sequences showing \emph{random} behavior; this is justified by the Glivenko--Cantelli theorem, which asserts that for a sequence $(U_k)_{k \geq 1}$ of independent, uniformly $[0,1]$-distributed random variables we have
$$
\lim_{N \to \infty} \frac{1}{N} \sum_{k=1}^N \mathbf{1}_{[a,b)} (U_k) = b-a \qquad \textrm{almost surely}
$$
for all $[a,b) \subset [0,1]$. Thus a deterministic sequence $(x_k)_{k \geq 1}$ which is u.d. mod 1 can be seen as a \emph{typical} realization of a random (uniformly $[0,1]$-distributed) sequence.\\

The theory a uniform distribution was boosted by Weyl's \cite{weyl} seminal paper of 1916, which contains the celebrated \emph{Weyl criterion} for uniform distribution of a sequence: a sequence$(x_k)_{k \geq 1}$ is u.d. mod 1 if and only if for all integers $h \neq 0$
\begin{equation} \label{weylk}
\lim_{N \to \infty} \frac{1}{N} \sum_{k=1}^N e^{2 \pi i h x_k} = 0.
\end{equation}
This criterion can be used to give an easy proof for the fact that the sequence $(\langle k x \rangle)_{k \geq 1}$ is u.d. mod 1 if and only if $x \not\in \mathbb{Q}$ (here $\langle \cdot \rangle$ stands for the fractional part function; usually this sequence is called $n \alpha$ rather than $k x$, but for the sake of consistency of the notation with later parts of this article we will denote it by $k x$). In fact, assume that $x \not\in \mathbb{Q}$; then, using the well-known formula for the geometric series, we have
$$
\frac{1}{N} \sum_{k=1}^N e^{2 \pi i h k x} = \frac{1}{N} ~ \frac{e^{2 \pi i h N x} - e^{2 \pi i h x}}{e^{2 \pi i h x} - 1} \to 0 \quad \textrm{as $N \to \infty$},
$$
where we used the fact that $e^{2 \pi i h x} - 1 \neq 0$ for $h \neq 0$ and $x \not\in \mathbb{Q}$. It is easy to see that for $x \in \mathbb{Q}$ the sequence $(\langle k x \rangle)_{k \geq 1}$ is \emph{not} u.d. modulo 1; thus the problem of deciding for which $x \in [0,1]$ the parametric sequence $(k x)_{k \geq 1}$ is u.d. mod 1 is completely solved.\\

Weyl's paper also contained a general result for parametric sequences of the form $(\langle n_k x\rangle)_{k \geq 1}$, where $(n_k)_{k \geq 1}$ is a sequence of distinct positive integers and $x$ is a real number from $[0,1]$: for almost all $x$ (in the sense of Lebesgue measure) the sequence $(\langle n_k x \rangle)_{k \geq 1}$ is u.d. mod 1. Accordingly the general case $(\langle n_k x \rangle)_{k \geq 1}$ resembles the properties of the case $(\langle k x \rangle)_{k \geq 1}$ insofar as in both cases the exceptional set is of measure zero; however, while in the latter case the exceptional set can be explicitly determined, it is generally very difficult to decide whether for a given sequence $(n_k)_{k \geq 1}$ and a given parameter $x$ the sequence $(\langle n_k x \rangle)_{k \geq 1}$ is u.d. mod 1 or not (see also Section \ref{sec2}).\\

In the following paragraph, we want to prove Weyl's result that $(\langle n_k x\rangle)_{k \geq 1}$ is u.d. mod 1 for almost all $x$. Throughout this article, we will repeatedly use methods from probability theory; this makes perfect sense, since the unit interval, equipped with Borel sets and Lebesgue measure, is a probability space (that is, a measure space $(\Omega, \mathcal{F}, \mathbb{P})$ for which $\mathbb{P}(\Omega) = 1$). 
We will use the Rademacher--Menshov inequality, which states that for a real orthonormal system $\phi_1(x), \dots, \phi_N(x)$ and for real coefficients $\alpha_1, \dots, \alpha_N$ we have
\begin{equation} \label{rm}
\int_0^1 \max_{1 \leq M \leq N} \left( \sum_{k=1}^M \alpha_k \phi_k \right)^2 ~dx \leq (\log_2 N + 2)^2 \sum_{k=1}^N \alpha_k^2
\end{equation}
(this inequality has been obtained independently by Rademacher \cite{rade} and Menshov \cite{mensh}; it can be proved quite easily using a dyadic splitting method, see e.g. \cite{loeve}). Note that an equivalent formulation of the Weyl criterion \eqref{weylk} is
$$
\lim_{N \to \infty} \frac{1}{N} \sum_{k=1}^N \cos 2 \pi h x_k = 0 \qquad \textrm{and} \qquad \lim_{N \to \infty} \frac{1}{N} \sum_{k=1}^N \sin 2 \pi h x_k = 0
$$
for all integers $h \neq 0$. For integers $m \geq 1$ and $h \neq 0$ we set
$$
S_{m,h} = \left\{ x \in [0,1]:~\max_{1 \leq M \leq 2^m} \left| \sum_{k=1}^M \cos 2 \pi h n_k x \right| > 2^{m/2} m^{2} \right\}.
$$
Then by \eqref{rm}, by the orthogonality of the trigonometric system and by the fact that by assumption the numbers $(n_k)_{k \geq 1}$ are distinct we have
\begin{equation} \label{che}
\int_0^1 \max_{1 \leq M \leq 2^m} \left( \sum_{k=1}^M \cos 2 \pi h n_k x \right)^2 dx \ll 2^m (\log 2^m)^2 \ll 2^m m^2
\end{equation}
(where ``$\ll $'' is the Vinogradov symbol). Chebyshev's inequality states that for any square-integrable function $f$ on $[0,1]$ we have that for any $t > 0$ 
\begin{equation} \label{cheby}
\lambda \Big(x \in [0,1]:~|f(x)| \geq t  \Big) \leq \frac{1}{t^2} \int_0^1 f(x)^2~dx.
\end{equation}
Applying this inequality, by \eqref{che} we have
$$
\lambda(S_{m,h}) \ll \frac{m^2}{m^4} \ll \frac{1}{m^2}.
$$
Thus
$$
\sum_{m=1}^\infty \lambda(S_{m,h}) < \infty,
$$
which by the first Borel--Cantelli lemma implies that with probability one (with respect to the Lebesgue measure on $[0,1]$) only finitely many events $S_{m,h}$ occur. Thus for almost all $x \in [0,1]$ there exists an $m_0 = m_0(x)$ such that
$$
\max_{1 \leq M \leq 2^m} \left| \sum_{k=1}^M \cos 2 \pi h n_k x \right| \leq 2^{m/2} m^{2} 
$$
for all $m \geq m_0$; consequently, there also exists an $N_0=N_0(x)$ such that
$$
\left| \sum_{k=1}^N \cos 2 \pi h n_k x \right| \leq 2 N^{1/2} (\log_2 N)^{2} 
$$
for $N \geq N_0$, which implies
$$
\lim_{N \to \infty} \frac{1}{N} \sum_{k=1}^N \cos 2 \pi h n_k x = 0.
$$
The same argument applies if we replace the function $\cos$ by $\sin$. Consequently, $(\langle n_k x \rangle)_{k \geq 1}$ is u.d. mod 1 for almost all $x \in [0,1]$.

\section{Metric number theory} \label{sec2}

One of the aims of metric number theory is to describe properties which are \emph{typical} for (real) numbers, where ``typical'' means that the exceptional set of numbers not possessing this property is small; in our case, we consider a property to be ``typical'' if it holds for almost all numbers in the sense of Lebesgue measure (but of course there are also other possibilities of deciding what a ``typical'' property is, for example by means of the Hausdorff dimension).\\

An early result from metric number theory is due to Borel: he proved that almost all numbers are \emph{normal} with respect to a given base $b \geq 2$ ($b$ being an integer). Here a real number $x \in [0,1]$ is called ``normal'' if in its base-$b$ expansion
$$
x = \sum_{i=1}^\infty r_i b^{-i}
$$
each digit $0, 1, \dots, b-1$ appears asymptotically with frequency $b^{-1}$, each block of~2 digits appears asymptotically with frequency $b^{-2}$, and, generally, each block of $d$ digits appears with asymptotic frequency $b^{-d}$. Formally this can be written as
\begin{equation} \label{beta}
\lim_{N \to \infty} \frac{1}{N} \sum_{k=1}^N \mathbf{1}_{[a b^{-d},(a+1)b^{-d})} (\langle b^{k-1} x \rangle) = b^{-d} 
\end{equation}
for all integers $d \geq 1$ and all integers $a \in \{0, \dots, b^d -1\}$. Historically, Borel's result is the first appearance of what we call today the \emph{strong law of large numbers}. To see this, we choose for simplicity $b=2$, and let for $x \in [0,1]$ the function $r_k(x)$ be defined as the $k$-th digit (after the decimal point) of the binary expansion of $x$. Then it an easy exercise to check that the functions $(r_k(x))_{k \geq 1}$, interpreted as random variables over the probability space $([0,1], \mathcal{B}([0,1]),\lambda)$, form a sequence of \emph{independent, identically distributed (i.i.d.) random variables} (remember that by definition a random variable is just a measurable function). Thus Borel's theorem, which in the special case $d=1$ (that is, for single digits) states that
$$
\lim_{N \to \infty} \frac{1}{N} \sum_{k=1}^N r_k(x) = \frac{1}{2} = \int_0^1 r_1(x) ~dx = \mathbb{E}_\lambda (r_1) \qquad \textrm{a.e.}
$$
is the strong law of large numbers for i.i.d. fair Bernoulli-distributed random variables (the functions $(r_k(x))_{k \geq 1}$ are called \emph{Rademacher functions}; the fact that they form a system if i.i.d. random variables has been first observed by Steinhaus in the 1920s). We just note by the way that Borel's theorem can also be interpreted as an early appearance of the pointwise ergodic theorem (for the transformation $Tx = \langle b x \rangle$).\\

Written in the form \eqref{beta} (which is not Borel's original notation) it is quite obvious that there is a connection between normal numbers and the criterion for uniform distribution modulo 1 in equation \eqref{ud}. Surprisingly, this connection was not noted (or, at least, not rigorously proved) before 1949, when Wall \cite{wall} showed that a number $x$ is normal in a base $b$ if and only if the sequence $(\langle b^k x \rangle)_{k \geq 1}$ is u.d. mod 1. Thus Borel's theorem can be also seen as a special case of Weyl's metric theorem on the uniform distribution of $( \langle n_k x \rangle)_{k \geq 1}$ for a.e. $x$.\\

Now we know that almost all numbers are normal (which was not so difficult to establish); on the other hand, constructing a normal is rather difficult, and checking whether a given number is normal or not is (usually) absolutely infeasible. Most constructions of normal numbers are based on the principle of concatenating blocks of digits generated by ``simple'' functions; for example, Champernowne's number (in base 10)
$$
0. ~1 ~2~3~4~5~6~7~8~9~10~11~12~\dots
$$
is obtained by concatenating the decimal expansions of the positive integers in consecutive order, the Copeland--Erd\H os number (again in base 10) is obtained by concatenating the decimal expansions of the primes
$$
0.~2~3~5~7~11~13~17~19~23~29~\dots,
$$
and there are several other constructions of this type (for example concatenating the values of polynomials \cite{nak} or other entire functions \cite{mad}). As mentioned before, checking whether a given number is normal or not is extremely difficult, and it is unknown whether constants such as $\sqrt{2}, e, \pi$ are normal. It is conjectured that all algebraic irrationals are normal, but no example or counterexample is known. For more details on this problem, see \cite{bc}. A closely related problem concerns sequences of the form $(\langle x^k \rangle)_{k \geq 1}$. By a result of Koksma \cite{koksma} this sequence is u.d. mod 1 for almost all $x>1$; however, not a single explicit value of $x$ for which this is the case is known. The sequence $(\langle (3/2)^k \rangle)_{k \geq 1}$ has attracted particular attention, but it is not even known whether $\limsup_{k \to \infty} \langle (3/2)^k \rangle - \liminf_{k \to \infty} \langle (3/2)^k \rangle \geq 1/2$ (\emph{Vijayaraghavan's conjecture} of 1940). For more information concerning this problem see \cite{dub1,dub2}.\\

Probably the most important open problem in metric number theory is the \emph{Duffin--Schaeffer conjecture} in metric Diophantine approximation. For a non-negative function $\psi: ~ \mathbb{N} \mapsto \mathbb{R}$, let $W(\psi)$ denote the set of real numbers $x \in [0,1]$ for which the inequality $|n x - a| < \psi(n)$  has infinitely many coprime solutions $(a,n)$. It is an easy application of the first Borel--Cantelli lemma to prove that $\lambda(W(\psi))=0$ if 
\begin{equation} \label{conv}
\sum_{n=1}^\infty \frac{\psi(n) \varphi(n)}{n} < \infty
\end{equation}
(here $\varphi$ denotes the Euler totient function); that means, divergence of the sum in \eqref{conv} is a \emph{necessary} condition for $\lambda(W(\psi))=1$. The Duffin--Schaeffer conjecture, proposed by R.J. Duffin and A.C. Schaeffer \cite{ds} in 1941, asserts that divergence of the sum in \eqref{conv} is also \emph{sufficient} to have $\lambda(W(\psi))=1$. Several special cases of the conjecture have been established (see for example \cite{harman}), but a complete solution of the problem seems to be far out of reach.\\

More information on the problems discussed in this section can be found in the books of Bugeaud \cite{bugeaud} and Harman \cite{harb}.

\section{Discrepancy}

The notion of the \emph{discrepancy} of a sequence has been introduced as a measure of the \emph{quality} of the uniform distribution mod 1 of a sequence. For a finite sequence $(x_1, \dots, x_N)$ of points in the unit interval, the discrepancy $D_N$ and the star-discrepancy $D_N^*$ are defined as
$$
D_N (x_1, \dots, x_N) = \sup_{0 \leq a < b \leq 1} \left|  \frac{\sum_{k=1}^N \mathbf{1}_{[a,b)} (x_k)}{N} - (b-a) \right|
$$
and
$$
D_N^* (x_1, \dots, x_N) = \sup_{0 \leq a \leq 1} \left|  \frac{\sum_{k=1}^N \mathbf{1}_{[0,a)} (x_k)}{N}  - a \right|.
$$
It is easy to see that these two discrepancies are equivalent in the sense that always $D_N^* \leq D_N \leq 2 D_N^*$, and that an infinite sequence $(x_k)_{k \geq 1}$ is u.d. mod 1 if and only if $D_N^*(x_1, \dots, x_N) \to 0$ as $N \to \infty$. An important inequality to estimate the discrepancy of a sequence is the \emph{Erd\H os--Tur\'an inequality}, which (in one out of many possible formulations) states that for any positive integer $H$
\begin{equation} \label{et}
D_N^*(x_1, \dots, x_N) \leq \frac{3}{H} + 3 \sum_{h=1}^H \frac{1}{h} \left| \sum_{k=1}^N e^{2 \pi i h x_k} \right|.
\end{equation}
We will use this inequality in Section \ref{sec5} to obtain an upper bound for the discrepancy of $(\langle n_k x \rangle)_{k \geq 1}$ for almost all $x$, by this means establishing a quantitative version of the theorem of Weyl mentioned in Section \ref{sec1}. Another important inequality concerning the discrepancy of sequences of points is \emph{Koksma's inequality}, which states that for any function $f$ which has bounded variation $\textup{Var}(f)$ in the unit interval the estimate
\begin{equation} \label{koks}
\left| \frac{1}{N} \sum_{k=1}^N f(x_k) - \int_0^1 f(x)~dx \right| \leq \textup{Var}(f) \cdot D_N^*(x_1, \dots, x_N)
\end{equation}
holds. The notions of uniform distribution mod 1 and discrepancy can be generalized in a natural way to the multi-dimensional setting, as can be the Erd\H os--Tur\'an inequality and Koksma's inequality (then called the \emph{Koksma-Hlawka inequality}). The multi-dimensional version of \eqref{koks} forms the foundation of the so-called \emph{Quasi-Monte Carlo method}, which is based on the observation that sequences having small discrepancy can be used for numerical integration.\\

By a result of Schmidt \cite{schm} there exists a positive constant $c$ such that for any infinite sequence $(x_k)_{k \geq 1}$ of points in the unit interval the inequality
$$
D_N^*(x_1, \dots, x_N) > c \frac{\log N}{N}
$$
holds for infinitely many $N$; on the other hand, there exist several constructions of sequences satisfying $D_N^*(x_1, \dots, x_N) = \mathcal{O}((\log N) N^{-1})$ as $N \to \infty$, so in the one-dimensional case the problem of the optimal asymptotic order of the discrepancy is solved. On the contrary, determining the optimal asymptotic order of the discrepancy in the multi-dimensional case turned out to be a very difficult problem, which is still open (see \cite{bi} for a survey).\\

More information on discrepancy theory and the Quasi-Monte Carlo method can be found in the books of Dick and Pillichshammer \cite{dpd} and Drmota and Tichy \cite{dts}. 

\section{Lacunary series}

The word \emph{lacunary} originates from the Latin \emph{lacuna} (ditch, gap), which is the diminutive form of \emph{lacus} (lake). Accordingly, a lacunary Fourier series is a series which has ``gaps'' in the sense that it is composed of trigonometric functions whose frequencies are far apart from each other. A classical gap condition is the \emph{Hadamard gap condition}, requiring that
\begin{equation} \label{had}
\frac{n_{k+1}}{n_k} \geq q > 1, \qquad k \geq 1;
\end{equation}
thus a (Hadamard) lacunary Fourier series is of the form
\begin{equation} \label{ko}
\sum_{k=1}^\infty \Big( a_k \cos 2 \pi n_k x + b_k \sin 2 \pi n_k x \Big)
\end{equation}
for $(n_k)_{k \geq 1}$ satisfying \eqref{had}. By a classical heuristics, lacunary sequences resemble many properties which are typical for sequences of \emph{independent} random variables. For example, by Kolmogorov's three series theorem a sequence of centered and uniformly bounded independent random variables $(X_k)_{k \geq 1}$ is almost surely convergent if and only if the variances satisfy
\begin{equation} \label{kol}
\sum_{k=1}^\infty \mathbb{V} (X_k) < \infty,
\end{equation}
and by a counterpart for lacunary series, also due to Kolmogorov, the series \eqref{ko} is almost everywhere convergent if and only if
\begin{equation} \label{carl}
\sum_{k=1}^\infty (a_k^2+b_k^2) < \infty. 
\end{equation}
Note here that the variance of the function $a_k \cos 2 \pi n_k x + b_k \sin 2 \pi n_k x$, considered as a random variable over the probability space $([0,1], \mathcal{B}([0,1]),\lambda)$, is simply given by
$$
\int_0^1 \left( a_k \cos 2 \pi n_k x + b_k \sin 2 \pi n_k x \right)^2 dx = a_k^2 + b_k^2.
$$
Thus the a.s. convergence behavior of series of independent random variables and of lacunary trigonometric series are in perfect accordance. Many similar results of the same type exist: for example, by a classical result of Salem and Zygmund \cite{sz}, under the gap condition \eqref{had} we have
$$
\lambda \left(x \in [0,1]:~\sum_{k=1}^N \cos 2 \pi n_k x < t \sqrt{N/2} \right) \to \Phi(t),
$$
where $\Phi(t)$ denotes the standard normal distribution function. In other words, the system $(\cos 2 \pi n_k x)_{k \geq 1}$ satisfies the central limit theorem. By a result of Erd\H os and G\'al \cite{egot}, the same system also satisfies the law of the iterated logarithm (LIL), that is
$$
\limsup_{N \to \infty} \frac{\left| \sum_{k=1}^N \cos 2 \pi n_k x\right|}{\sqrt{2 N \log \log N}} = \frac{1}{\sqrt{2}} \qquad \textup{a.e.}
$$

The situation gets significantly more complicated if we consider the more general sequence $(f(n_k x))_{k \geq 1}$ for a (in some sense) ``nice'' function $f$ satisfying
\begin{equation} \label{f}
f(x+1) = f(x), \qquad \int_0^1 f(x) ~dx=0,
\end{equation}
instead of $(\cos 2 \pi n_k x)_{k \geq 1}$. A striking result for this general setting is a theorem of Philipp \cite{ph}, who confirmed the so-called Erd\H os--G\'al conjecture by proving that under \eqref{had} we have
\begin{equation} \label{phi}
\frac{1}{4 \sqrt{2}} \leq \limsup_{N \to \infty} \frac{N D_N^*(\{n_1 x\}, \dots, \{n_N x\})}{\sqrt{2 N \log \log N}} \leq C_q \qquad \textup{a.e.};
\end{equation}
this is a counterpart of the Chung--Smirnov LIL for the Kolmogorov--Smirnov statistic in probability theory. As a consequence of \eqref{phi} and Koksma's inequality \eqref{koks} we have
\begin{equation} \label{phi2}
\limsup_{N \to \infty} \frac{ \left| \sum_{k=1}^N f(n_k x) \right|}{\sqrt{2 N \log \log N}}  \leq C_{f,q} \qquad \textup{a.e.}
\end{equation}
Calculating the precise value of the $\limsup$ in \eqref{phi} and \eqref{phi2} is a very difficult problem, and depends on number-theoretic properties of $(n_k)_{k \geq 1}$ and Fourier-analytic properties of $f$ (or of the indicator functions in the case of \eqref{phi}) in a very delicate way. In the case of $n_k = \theta^k$ for an integer $\theta$ the problem has been solved by Fukuyama \cite{fuku}; he proved that almost everywhere 
$$
\limsup_{N \to \infty} \frac{N D_N^*(\{\theta x\}, \dots, \{\theta^N x\})}{\sqrt{2 N \log \log N}} = \left\{ \begin{array}{ll} \frac{\sqrt{42}}{9} & \textrm{if $\theta=2$}\\ \frac{\sqrt{(\theta+1)\theta(\theta-2)}}{2 \sqrt{(\theta-1)^3}} & \textrm{if $\theta \geq 4$ is even} \\ \frac{\sqrt{\theta+1}}{2 \sqrt{\theta-1}} & \textrm{if $\theta \geq 3$ is odd} \end{array} \right.
$$
In view of the results mentioned in Section \ref{sec2}, Fukuyama's theorem establishes the typical asymptotic order of the discrepancy of normal numbers. In a sense the sequence $(\theta^k)_{k \geq 1}$ is a pathological example of a lacunary sequence, exhibiting an extremely strong relation between its consecutive terms. For a lacunary sequence for which no such strong arithmetic relations exist, the LIL is satisfied in the form 
$$
\limsup_{N \to \infty} \frac{N D_N^*(\{n_1 x\}, \dots, \{n_N x\})}{\sqrt{2 N \log \log N}} = \frac{1}{2} \qquad \textrm{a.e.},
$$
which is in perfect accordance (including the value of the constant of the right-hand side) with the Chung--Smirnov LIL for i.i.d. random variables (see \cite{alil} for details). For more information on lacunary sequences in the context of metric discrepancy theory and probabilistic limit theorems, see the survey paper \cite{absurv}.\\

Lacunary functions are well-known in analysis for several other interesting properties, apart from their resemblance of the behavior of systems of independent random variables. For example, Weierstrass's celebrated example of a nowhere differentiable function is defined by means of a lacunary trigonometric series. It should be noted that the notion of lacunary series does not only include lacunary trigonometric series, but also other series such as for example lacunary Taylor series. For a survey, see \cite{kah}.

\section{Almost everywhere convergence} \label{sec5}

The Kolmogorov three series theorem gives a full characterization of the a.s. convergence behavior of sums of independent random variables. In general the a.s. convergence condition comprises of three conditions about the convergence or divergence of certain series (which explains the name \emph{three series theorem}), but in the case of centered, uniformly bounded random variables the criterion reduces to the simple condition \eqref{kol}.\\

As noted in the previous section, there exists an analogue of the three series theorem for the case of Hadamard lacunary trigonometric series; however, surprisingly, the requirement of considering a Fourier series which contains only frequencies along an exponentially growing subsequence can be entirely dropped. This is the celebrated Carleson's theorem \cite{carl}, which is considered as one of the major achievements of Fourier analysis in 20th-century mathematics: a Fourier series
$$
\sum_{k=1}^\infty \Big( a_k \cos 2 \pi k x + b_k \sin 2 \pi k x \Big)
$$
is a.e. convergent, provided \eqref{carl} holds. Moreover, for any function $f \in L^2([0,1])$ and
$$
f(x)  \sim \sum_{k=1}^\infty \Big( a_k \cos 2 \pi k x + b_k \sin 2 \pi k x \Big),
$$
setting
$$
s_N(f;x) = \sum_{k=1}^N \Big( a_k \cos 2 \pi k x + b_k \sin 2 \pi k x \Big)
$$
we have
$$
s_N(f;x) \to f(x) \qquad \textrm{as $N \to \infty$} \qquad \textrm{for a.e. $x$}.
$$
Carleson's theorem has a breathtaking consequence: there exists an absolute constant $c$ such that for any function $f$ in $L^2([0,1])$ we have, writing $\| \cdot \|$ for the $L^2([0,1])$ norm,
\begin{equation} \label{ch}
\left\| \sup_{N \geq 1} \left| s_N(f;x) \right| \right\| \leq c \|f\|.
\end{equation}
Carleson's theorem has been extended to the case $f \in L^p([0,1]),~p>1,$ by Hunt \cite{hunt}. It is a very deep result, and although alternative proofs have been given by Fefferman \cite{feff} and Lacey and Thiele \cite{lt}, no ``easy'' proof exists. For a comprehensive treatment of the subject, see the monograph of Arias de Reyna \cite{arias} and the survey paper \cite{lac}.\\

As an application of Carleson's theorem, we will show how it can be used to obtain a quantitative version of the results on a.e. uniform distribution of Weyl mentioned in Section \ref{sec1}. This argument is due to Baker \cite{baker} and leads to the upper bound
\begin{equation} \label{baker}
D_N^*(\{n_1 x\}, \dots, \{n_N x\}) \ll \frac{(\log N)^{3/2+\varepsilon}}{\sqrt{N}} \qquad \textup{a.e.}
\end{equation}
for any strictly increasing sequence of positive integers $(n_k)_{k \geq 1}$ and any $\varepsilon>0$. To outline the similarity between this proof and the one given in Section \ref{sec1}, we will use a real version of the Erd\H os--Tur\'an inequality \eqref{et}.\\

Let $\varepsilon>0$. For integers $m \geq 1$ we set
$$
S_{m} = \left\{ x \in [0,1]:~\max_{1 \leq M \leq 2^m} \left|\sum_{1 \leq h \leq 2^{m/2}} \frac{1}{h} \sum_{k=1}^M \cos 2 \pi h n_k x \right| > 2^{m/2} m^{3/2+\varepsilon} \right\}.
$$
Note that by \eqref{ch} for any $h$ we have
\begin{equation} \label{cha}
\left\| \max_{1 \leq M \leq 2^m} \left| \sum_{k=1}^M \cos 2 \pi h n_k x \right| \right\| \leq c \left\| \sum_{k=1}^{2^m} \cos 2 \pi h n_k x \right\| \leq c 2^{m/2}
\end{equation}
for an absolute constant $c$. Thus by Minkowski's inequality we have
\begin{eqnarray*}
& & \left\| \max_{1 \leq M \leq 2^m} \left|\sum_{1 \leq h \leq 2^{m/2}} \frac{1}{h} \sum_{k=1}^M \cos 2 \pi h n_k x \right| \right\| \\
& \leq & \sum_{1 \leq h \leq 2^{m/2}} \frac{1}{h}  \left\| \max_{1 \leq M \leq 2^m} \left|\sum_{k=1}^M \cos 2 \pi h n_k x \right| \right\| \\
& \ll & c m 2^{m/2}.
\end{eqnarray*}
Consequently, by Chebyshev's inequality \eqref{cheby} we have
$$
\lambda(S_{m}) \ll \frac{1}{m^{1+2\varepsilon}}.
$$
Thus
$$
\sum_{m=1}^\infty \lambda(S_{m}) < \infty,
$$
which by the first Borel--Cantelli lemma means that with probability one only finitely many events $S_{m}$ occur. Thus for almost all $x \in [0,1]$ there exists an $m_0 = m_0(x)$ such that
$$
\max_{1 \leq M \leq 2^m} \left|\sum_{1 \leq h \leq 2^{m/2}} \frac{1}{h} \sum_{k=1}^M \cos 2 \pi h n_k x \right| \leq 2^{m/2} m^{3/2+\varepsilon} 
$$
for all $m \geq m_0$; consequently, there also exists an $N_0=N_0(x)$ such that
$$
\left| \sum_{1 \leq h \leq \sqrt{N}} \frac{1}{h} \sum_{k=1}^N \cos 2 \pi h n_k x \right|  \leq 2 N^{1/2} (\log_2 N)^{3/2+\varepsilon} 
$$
for $N \geq N_0$. The same result holds if we replace the function $\cos$ by $\sin$, and using \eqref{et} (split into real and imaginary part) we get \eqref{baker}. Carleson's inequality \eqref{ch} in the form \eqref{cha} plays a key role in this proof, and if it is replaced by the Rademacher--Menshov inequality, which gives an additional logarithmic factor as in \eqref{che}, one can only obtain \eqref{baker} with the exponent $3/2+\varepsilon$ replaced by $5/2+\varepsilon$.\\

The optimal exponent of the logarithmic term in \eqref{baker} is an important open problem in metric discrepancy theory. Note that by Koksma's inequality \eqref{koks} as a consequence of \eqref{baker} we get
\begin{equation} \label{fvar}
\left| \sum_{k=1}^N f(n_k x) \right| \ll \sqrt{N} (\log N)^{3/2+\varepsilon} \qquad \textup{a.e.}
\end{equation}
for any function $f$ satisfying \eqref{f} which has bounded variation on $[0,1]$. On the other hand, Berkes and Philipp \cite{berkes} constructed a sequence $(n_k)_{k \geq 1}$ for which
$$
\limsup_{N \to \infty} \frac{\left| \sum_{k=1}^N \cos 2 \pi n_k x \right|}{\sqrt{N \log N}} = \infty \qquad \textup{a.e.}, 
$$
which again by Koksma's inequality implies that the exponent of the logarithmic term in \eqref{baker} can in general not be reduced below $1/2$. In the following section we will see that $1/2$ is in fact the optimal exponent, at least if we consider only a single function $f$ (as in \eqref{fvar}) and not the discrepancy $D_N^*$.\\

For more information on the a.e. convergence of sums of dilated functions, the interested reader is referred to the comprehensive survey article of Berkes and Weber~\cite{bewe}.

\section{Sums involving greatest common divisors}

In the context of counting lattice points in right-angled triangles, around 1920 Hardy and Littlewood investigated the problem of finding good upper bounds for the asymptotic order of
\begin{equation} \label{hl}
\sum_{k=1}^N  \left( \{k x \} -1/2 \right) \qquad \textrm{as} \qquad N \to \infty.
\end{equation}
The $L^2([0,1])$-norm of~\eqref{hl} can be calculated using the formula
\begin{equation} \label{land}
\int_0^1\left( \{m x \} -1/2 \right) \left( \{n x \} -1/2 \right) ~dx = \frac{1}{12} \frac{(\gcd(m,n))^2}{mn}
\end{equation}
for integers $m,n$ (first stated by Franel and proved by Landau in 1924). The generalized problem of estimating
\begin{equation} \label{galsum}
\sum_{k,l=1}^N \frac{(\gcd(n_k,n_l))^2}{n_k n_l} 
\end{equation}
for an arbitrary sequence of distinct positive integers $n_1, \dots, n_N$ was posed as a prize problem by the Scientific Society at Amsterdam in 1947 (following a suggestion of Erd\H os), and solved by G\'al~\cite{gal} in 1949. He proved that there exists an absolute constant $c$ such that~\eqref{galsum} is bounded by $c N (\log \log N)^2$, and that this upper bound is asymptotically optimal. Koksma~\cite{koks} observed that as a consequence for any centered, periodically extended indicator function $f(x)=\mathbf{1}_{(a,b)} (\{x\}) - (b-a)$ (and in fact even for any 1-periodic function $f$ having mean zero and bounded variation on [0,1]) the estimate
\begin{equation} \label{koko}
\int_0^1 \left(\sum_{k=1}^N f(n_k x)\right)^2~dx \ll N (\log \log N)^2
\end{equation}
holds. This follows from a generalized version of \eqref{land}, which we will deduce in the next few lines. Assume that $f$ satisfies \eqref{f} and is of bounded variation on $[0,1]$. Let
$$
f(x) \sim \sum_{j=1}^\infty a_j \cos 2 \pi j x
$$
denote the Fourier series of $f$ (for simplicity we assume that it is a pure cosine-series; the general case works in exactly the same way). Then 
\begin{equation} \label{aj}
|a_j| \ll j^{-1}
\end{equation}
(see \cite[p.~48]{zt}; this estimate can be easily proved using the fact that any function of bounded variation can be written as the sum of two bounded and monotone functions). Thus for integers $m,n$ we have, by \eqref{aj} and the orthogonality of the trigonometric system, and writing $\delta(\cdot,\cdot)$ for the Kronecker function, 
\begin{eqnarray}
\int_0^1 f(mx) f(nx) ~dx & = & \sum_{j_1, j_2 = 1}^\infty a_{j_1} a_{j_2} \delta(m j_1, n j_2) \nonumber\\
& \ll & \sum_{j_1, j_2 = 1}^\infty \frac{1}{j_1 j_2} \delta(m j_1, n j_2). \label{sound}
\end{eqnarray}
Now $m j_1 = n j_2$ is only possible if $j_1 = j n/\gcd(m,n)$ and $j_2 = j m /\gcd(m,n)$ for some integer $j \geq 1$. Consequently, \eqref{sound} is at most
$$
\sum_{j=1}^\infty \frac{\gcd(m,n)}{jn} \frac{\gcd(m,n)}{jm} \ll \frac{\gcd(m,n)^2}{mn},
$$
which implies, together with the aforementioned result of G\'al, that \eqref{koko} holds. Sums involving common divisors similar to~\eqref{galsum} were studied by Dyer and Harman~\cite{dh} in the context of metric Diophantine approximation. They investigated 
\begin{equation} \label{gcda}
\max_{n_1 < \dots < n_N} \sum_{k,l=1}^N \frac{(\gcd(n_k,n_l))^{2 \alpha}}{(n_k n_l)^\alpha},  \qquad \alpha \in [1/2,1),
\end{equation}
and, amongst other results, proved for the particularly interesting case $\alpha=1/2$ the upper bound
\begin{equation} \label{gcda2}
\max_{n_1 < \dots < n_N} \sum_{k,l=1}^N \frac{\gcd(n_k,n_l)}{\sqrt{n_k n_l}} \ll N \exp \left( \frac{5 \log N}{\log \log N} \right).
\end{equation}
Recently, Aistleitner, Berkes and Seip~\cite{abs} obtained upper bounds for ~\eqref{gcda} which are essentially optimal. They proved that
\begin{equation}\label{abs}
\max_{n_1 < \dots < n_N} \sum_{k,l=1}^N \frac{(\gcd(n_k,n_l))^{2 \alpha}}{(n_k n_l)^\alpha}  \leq C_\varepsilon N \exp\left( (1+\varepsilon) g(\alpha,N) \right),
\end{equation}
where for $ 1/2 < \alpha < 1$ we have 
\begin{equation} \label{ga}
g(\alpha,N)=  \frac{\left(\frac{8}{1-\alpha}+\frac{16\cdot 2^{-\alpha}}{\sqrt{2\alpha-1}}\right) (\log N)^{1-\alpha}}{(\log\log N)^{\alpha}}
                                               + \frac{(\log N)^{(1-\alpha)/2}}{1-\alpha},
\end{equation} 
for $\alpha =1/2$ we have
\begin{equation} \label{ga2}
g(1/2,N) = 25 \sqrt{\log N}\sqrt{\log\log N},
\end{equation}
and where $C_\varepsilon$ is a constant only depending on $\varepsilon>0$. Here the asymptotic order of $g(\alpha,N)$ in \eqref{ga} is optimal, and the asymptotic order of $g(1/2,N)$ in \eqref{ga2} can perhaps be reduced from $\sqrt{\log N}\sqrt{\log \log N}$ to $\sqrt{\log N}/\sqrt{\log \log N}$ (but not below). As an application, Aistleitner, Berkes and Seip improved the exponent of the logarithmic term in~\eqref{fvar} to $1/2+\varepsilon$, which is optimal (up to the $\varepsilon$). Another application of such GCD sums is concerning the a.e. convergence of series $\sum_{k=1}^\infty c_k f(n_k x)$ for functions $f$ of bounded variation or being H\"older-continuous (see \cite{abs}; cf. also \cite{alip}) or of so-called Davenport series (see \cite{bre}). There is also a close connection with certain properties of the Riemann zeta function, which requires further investigation (cf. \cite{hilb}).\\

As a consequence of \eqref{ga} (and using a trick to modify the argument around equation \eqref{sound} in such a way to get a generalized GCD sum of the form \eqref{gcda} instead of \eqref{galsum}) one can show that for any $f$ of bounded variation satisfying \eqref{f} and for any strictly increasing sequence $(n_k)_{k \geq 1}$ of positive integers the Carleson-type inequality
$$
\left\| \max_{1 \leq M \leq N} \left| \sum_{k=1}^M c_k f(n_k x) \right| \right\| \leq c (\log \log N)^4 \sum_{k=1}^N c_k^2 
$$
holds, which, using an argument similar to the one used to prove \eqref{baker} in the previous section, leads to 
\begin{equation} \label{fvar2}
\left| \sum_{k=1}^N f(n_k x) \right| \ll \sqrt{N} (\log N)^{1/2+\varepsilon} \qquad \textup{a.e.}
\end{equation}
As mentioned at the end of the previous section this means that the problem concerning the optimal (a.e.) asymptotic order of the sum $\sum f(n_k x)$ is solved; however, the more difficult case of the precise a.e. asymptotic order of the discrepancy $D_N^*(\{n_1 x\}, \dots, \{n_N x\})$ remains open.\\

\emph{Concluding remark:} There exists a close connection between discrepancy theory and harmonic analysis, which we have for example observed in Weyl's criterion~\eqref{weylk}, in the Erd\H os--Tur\'an inequality~\eqref{et} and in the Carleson convergence theorem in Section~\ref{sec5}. This connection goes far beyond the material contained in this article, and is comprehensively presented in a survey article of Dmitriy Bilyk in the present volume (see~\cite{bilyksurv}).\footnote{This remark refers to the conference proceedings volume for the RICAM workshop on ``Uniform Distribution and Quasi-Monte Carlo Methods'' (held from October 14--18, 2013, in Linz, Austria), in which the present article will be published. The proceedings volume will appear as part of the ``Radon Series on Computational and Applied Mathematics'', published by DeGruyter.}


\begin{thebibliography}{10}

\bibitem{alil}
Christoph Aistleitner, On the law of the iterated logarithm for the discrepancy
  of lacunary sequences, \emph{Trans. Amer. Math. Soc.} {362} (2010),
  5967--5982.

\bibitem{alip}
\bysame, Convergence of {$\sum c_kf(kx)$} and the {L}ip {$\alpha$} class,
  \emph{Proc. Amer. Math. Soc.} {140} (2012), 3893--3903.

\bibitem{absurv} 
Christoph Aistleitner and Istvan Berkes, Probability and metric discrepancy theory,
 \emph{ Stoch. Dyn.} {11} (2011), 183--207.

\bibitem{abs}
Christoph Aistleitner, Istvan Berkes and Kristian Seip, {G}{C}{D} sums from
  {P}oisson integrals and systems of dilated functions, Preprint. Available at
  \url{http://arxiv.org/abs/1210.0741}.

\bibitem{arias}
Juan Arias~de Reyna, \emph{Pointwise convergence of {F}ourier series}, Lecture
  Notes in Mathematics 1785, Springer-Verlag, Berlin, 2002.

\bibitem{bc}
David~H. Bailey and Richard~E. Crandall, On the random character of fundamental
  constant expansions, \emph{Experiment. Math.} {10} (2001), 175--190.

\bibitem{baker}
R.~C. Baker, Metric number theory and the large sieve, \emph{J. London Math.
  Soc. (2)} {24} (1981), 34--40.

\bibitem{berkes}
Istv{\'a}n Berkes and Walter Philipp, The size of trigonometric and {W}alsh
  series and uniform distribution {${\rm mod}\ 1$}, \emph{J. London Math. Soc.
  (2)} {50} (1994), 454--464.

\bibitem{bewe}
Istv{\'a}n Berkes and Michel Weber, On the convergence of {$\sum c_kf(n_kx)$}, 
 \emph{Mem. Amer. Math. Soc.} {201} (2009).

\bibitem{bilyksurv}
Dmitriy Bilyk, Discrepancy theory and harmonic analysis, \emph{contained in this volume}.

\bibitem{bi}
\bysame, On {R}oth's orthogonal function method in discrepancy theory,
  \emph{Unif. Distrib. Theory} {6} (2011), 143--184.

\bibitem{bre}
Julien Br{\'e}mont, Davenport series and almost-sure convergence, \emph{Q. J.
  Math.} {62} (2011), 825--843.

\bibitem{bugeaud}
Yann Bugeaud, \emph{Distribution modulo one and {D}iophantine approximation},
  Cambridge Tracts in Mathematics 193, Cambridge University Press, Cambridge,
  2012.

\bibitem{carl}
Lennart Carleson, On convergence and growth of partial sums of {F}ourier
  series, \emph{Acta Math.} {116} (1966), 135--157.

\bibitem{dpd}
Josef Dick and Friedrich Pillichshammer, \emph{Digital nets and sequences},
  Cambridge University Press, Cambridge, 2010.

\bibitem{dts}
Michael Drmota and Robert~F. Tichy, \emph{Sequences, discrepancies and
  applications}, Lecture Notes in Mathematics 1651, Springer-Verlag, Berlin,
  1997.

\bibitem{dub1}
Art{\=u}ras Dubickas, On the distance from a rational power to the nearest
  integer, \emph{J. Number Theory} {117} (2006), 222--239.

\bibitem{dub2}
\bysame, On the powers of 3/2 and other rational numbers, \emph{Math. Nachr.}
  {281} (2008), 951--958.

\bibitem{ds} 
R.J. Duffin and A.C. Schaeffer, Khintchine's problem in metric Diophantine approximation,
 \emph{Duke Math. J.} {8} (1941), 243--255.

\bibitem{dh}
Tony Dyer and Glyn Harman, Sums involving common divisors, \emph{J. London
  Math. Soc. (2)} {34} (1986), 1--11.

\bibitem{egot}
P.~Erd{\"o}s and I.~S. G{\'a}l, On the law of the iterated logarithm. {I},
  {II}, \emph{Nederl. Akad. Wetensch. Proc. Ser. A. {\bf 58} = Indag. Math.}
  {17} (1955), 65--76, 77--84.

\bibitem{feff}
Charles Fefferman, Pointwise convergence of {F}ourier series, \emph{Ann. of
  Math. (2)} {98} (1973), 551--571.

\bibitem{fuku}
K.~Fukuyama, The law of the iterated logarithm for discrepancies of
  {$\{\theta^nx\}$}, \emph{Acta Math. Hungar.} {118} (2008), 155--170.

\bibitem{gal}
I.~S. G{\'a}l, A theorem concerning {D}iophantine approximations, \emph{Nieuw
  Arch. Wiskunde (2)} {23} (1949), 13--38.

\bibitem{harman}
Glyn Harman, Some cases of the {D}uffin and {S}chaeffer conjecture,
  \emph{Quart. J. Math. Oxford Ser. (2)} {41} (1990), 395--404.

\bibitem{harb}
\bysame, \emph{Metric number theory}, London Mathematical Society Monographs.
  New Series~18, The Clarendon Press Oxford University Press, New York, 1998.

\bibitem{hilb}
Titus Hilberdink, An arithmetical mapping and applications to
  {$\Omega$}-results for the {R}iemann zeta function, \emph{Acta Arith.} {139}
  (2009), 341--367.

\bibitem{hunt}
Richard~A. Hunt, \emph{On the convergence of {F}ourier series}, Orthogonal
  {E}xpansions and their {C}ontinuous {A}nalogues ({P}roc. {C}onf.,
  {E}dwardsville, {I}ll., 1967), Southern Illinois Univ. Press, Carbondale,
  Ill., 1968, pp.~235--255.

\bibitem{kah}
J.-P. Kahane, Lacunary {T}aylor and {F}ourier series, \emph{Bull. Amer. Math.
  Soc.} {70} (1964), 199--213.

\bibitem{koksma}
J.~F. Koksma, Ein mengentheoretischer {S}atz \"uber die {G}leichverteilung
  modulo {E}ins, \emph{Compositio Math.} {2} (1935), 250--258.

\bibitem{koks}
\bysame, On a certain integral in the theory of uniform distribution,
  \emph{Nederl. Akad. Wetensch., Proc. Ser. A. {\bf 54} = Indagationes Math.}
  {13} (1951), 285--287.

\bibitem{lac}
Michael Lacey, Carleson's theorem: proof, complements, variations, \emph{Publ.
  Mat.} {48} (2004), 251--307.

\bibitem{lt}
Michael Lacey and Christoph Thiele, A proof of boundedness of the {C}arleson
  operator, \emph{Math. Res. Lett.} {7} (2000), 361--370.

\bibitem{loeve}
Michel Lo{\`e}ve, \emph{Probability theory}, 2nd ed. The University Series in
  Higher Mathematics. D. Van Nostrand Co., Inc., Princeton, N. J.-Toronto-New
  York-London, 1960.

\bibitem{mad}
Manfred~G. Madritsch, J{\"o}rg~M. Thuswaldner and Robert~F. Tichy, Normality of
  numbers generated by the values of entire functions, \emph{J. Number Theory}
  {128} (2008), 1127--1145.

\bibitem{mensh}
D.~Menchoff, {Sur les s\'eries de fonctions orthogonales. (Premi\'ere Partie.
  La convergence.).}, \emph{Fundamenta math.} {4} (1923), 82--105 (French).

\bibitem{nak}
Yoshinobu Nakai and Iekata Shiokawa, Discrepancy estimates for a class of
  normal numbers, \emph{Acta Arith.} {62} (1992), 271--284.

\bibitem{ph}
Walter Philipp, Limit theorems for lacunary series and uniform distribution
  {${\rm mod}\ 1$}, \emph{Acta Arith.} {26} (1974/75), 241--251.

\bibitem{rade}
H.~Rademacher, {Einige S\"atze \"uber Reihen von allgemeinen
  Orthogonalfunktionen.}, \emph{Math. Ann.} {87} (1922), 112--138 (German).

\bibitem{sz}
R.~Salem and A.~Zygmund, On lacunary trigonometric series, \emph{Proc. Nat.
  Acad. Sci. U. S. A.} {33} (1947), 333--338.

\bibitem{schm}
Wolfgang~M. Schmidt, Irregularities of distribution. {VII}, \emph{Acta Arith.}
  {21} (1972), 45--50.

\bibitem{wall}
Donald~D. Wall, \emph{Normal numbers}, Ph.D. thesis, University of California,
  Berkeley, 1949.

\bibitem{weyl}
H.~Weyl, {\"Uber die Gleichverteilung von Zahlen mod. Eins.}, \emph{Math. Ann.}
  {77} (1916), 313--352 (German).

\bibitem{zt}
A.~Zygmund, \emph{Trigonometric series. {V}ol. {I}, {II}}, Cambridge
  Mathematical Library, Cambridge University Press, Cambridge, 1988, Reprint of
  the 1979 edition.

\end{thebibliography}

\providecommand{\bysame}{\leavevmode\hbox to3em{\hrulefill}\thinspace}
\providecommand{\MR}{\relax\ifhmode\unskip\space\fi MR }
\providecommand{\MRhref}[2]{%
  \href{http://www.ams.org/mathscinet-getitem?mr=#1}{#2}
}
\providecommand{\href}[2]{#2}

\end{document}